\documentclass[psamsfonts]{amsart}

\newtheorem{theorem}{Theorem}[section]
\newtheorem{lemma}[theorem]{Lemma}
\newtheorem{proposition}[theorem]{Proposition}
\newtheorem{corollary}[theorem]{Corollary}

\theoremstyle{definition}
\newtheorem{definition}[theorem]{Definition}
\newtheorem{example}[theorem]{Example}

\newtheorem{remark}[theorem]{Remark}

\newtheorem{problem}[theorem]{Problem}

\theoremstyle{remark}

\numberwithin{equation}{section}

\newcommand{\NN}{\mathbb{N}}
\newcommand{\ZZ}{\mathbb{Z}}
\newcommand{\QQ}{\mathbb{Q}}
\newcommand{\RR}{\mathbb{R}}
\newcommand{\CC}{\mathbb{C}}

\newcommand{\PP}{\mathbb{P}}
\renewcommand{\AA}{\mathbb{A}}

\newcommand  {\shL}     {\mathcal{L}}

\newcommand  {\foa}     {\mathfrak{a}}

\newcommand  {\fob}     {\mathfrak{b}}

\newcommand  {\fop}     {\mathfrak{p}}

\newcommand  {\foq}     {\mathfrak{q}}

\newcommand  {\for}     {\mathfrak{r}}

\newcommand  {\Cl}      {\operatorname{Cl}}

\newcommand  {\Char}    {\operatorname{char}}

\newcommand  {\Div}     {\operatorname{Div}}

\newcommand  {\Num}     {\operatorname{Num}}

\newcommand  {\primid}  {\mathfrak{p}}
\renewcommand{\O}       {\mathcal{O}}
\renewcommand{\OE}      {ohne Einschrankung}

\newcommand  {\Pic}     {\operatorname{Pic}}

\newcommand  {\Proj}    {\operatorname{Proj}}

\newcommand  {\ra}      {\rightarrow}

\newcommand  {\Spec}    {\operatorname{Spec}}

\newcommand  {\Supp}    {\operatorname{Supp}}

\def\mydate{\number\day\space\ifcase\month \or January\or February\or March\or April\or May\or
June\or July\or August\or September\or October\or November\or
December\fi \space\number\year}

\usepackage{amscd}
\usepackage{amssymb}

\begin{document}

\title[The affine Class group of a normal scheme]
{The affine Class group of a normal scheme}

% Remove or comment out any unused author tags.
% author one information

\author[Holger Brenner]{Holger Brenner}
\address{Mathematische Fakult\"at, Ruhr-Universit\"at, 
               44780 Bochum, Germany}
%\curraddr{}
\email{brenner@cobra.ruhr-uni-bochum.de}

%\thanks{}

\subjclass{}
%\date{}

% at present the "communicated by" line appears only in ERA and PROC
%\commby{}

%\dedicatory{Preliminary version,  \mydate}

\begin{abstract}
We study the property of a normal scheme, that the complement of every
hypersurface is an affine scheme.
To this end we introduce the affine class group. It is a factor group
of the divisor class group and measures the deviation from this property.
We study the behaviour of the affine class
group under faithfully flat extensions
and under the formation of products, and we compute it for different
classes of rings.
\end{abstract}

\maketitle

\noindent
Mathematical Subject Classification (1991):
13C20; 13F15; 14C20; 14C25; 32E10

%===========================================================
\section*{Introduction}
Consider a hypersurface $V \subseteq X=\Spec\, A$, where $A$ is a
factorial domain.
Then $V=V(f)$ and the open complement $X-V=D(f)\cong \Spec\, A_f$ is again
an affine scheme.
Thus the spectrum $X$ of a factorial domain has 
the property that the complement of every hypersurface is an affine
scheme. This property holds also in the following situations.
\par\smallskip\noindent
-$X$ is locally $\QQ$-factorial,
meaning that the divisor class group of each point in $X$ is a torsion group.
This rests upon the fact that affineness of a morphism
is a local property on the base.

\par\smallskip\noindent
-$X$ is regular. For then $X$ is locally factorial.

\par\smallskip\noindent
-$X$ is a curve, i.e. a one dimensional noetherian scheme. This is a
consequence of the Hauptidealsatz of Krull.

\par\smallskip\noindent
-$X$ is the spectrum of a normal excellent domain $A$ of dimension two.
The first proof of this result was given by Nagata in connection with the 14th
problem of Hilbert, see \cite{naghil} and \cite{binsto1}.
\par\smallskip

In this paper we shall study this property of a normal
noetherian affine scheme.
For this purpose we introduce in the first section the
\emph{affine class group}
of $A$, $ACl \, A$. It is
a torsion free residue class group of the divisor class group $Cl \, A$
and it measures the deviation of the
described property in a similar way as $Cl \, A$
measures the deviation from factoriality.

In the second section we study the behavior of the affine class group under
faithfully flat ring extensions.
We show that the affine class group of a normal domain $A$ and
of the polynomial ring $A[T]$ coincide (Cor. \ref{affinclassgauss}).
In particular, if our property holds for $A$, it holds for $A[T]$.
For a normal excellent domain $A$
of dimension two we obtain that
the complement of every (hyper)surface in $\Spec\, A[T]$ is affine
(Cor. \ref{nagatazylinder}). A geometric consequence of this result
is that in the affine line over a normal affine surface
the intersection of two surfaces cannot contain isolated points. 

In the third section we study hypersurfaces
in the product of a normal affine variety $X$ with an affine smooth
curve $C$ over an algebraically closed field.
The crucial point here is to look at the closure of the graph
of a rational mapping from $X$ to $C$.
For a normal surface $X$ we obtain again
that the complement of every hypersurface
in $X \times C$ is affine (Theorem \ref{surfacecurve}).

The aim of section four and five is to compute the affine class group
in some examples, in order to illustrate how $ACl\, A$ and $Cl\, A$
are related.
The affine class group of hyperbolas,
of monoid rings and of determinantal rings is just
the divisor class group modulo torsion.

Section five deals with relations between affineness properties of a
smooth projective variety and of an affine cone over it.
Under suitable conditions the pull back of a numerically trivial divisor
vanishes in the affine class group.
For example we get that the 
affine class group of an affine cone over a
geometrically ruled surface is always $\ZZ$, thereas the divisor class group
may be very big.

%===========================================================
\section{Definitions}
Throughout this paper, by a normal scheme we shall mean a
noetherian separated irreducible reduced normal scheme.
A divisor is always a Weil divisor.
The complement of an effective divisor is the complement
of its support.

\begin{definition}
Let $X$ denote a normal scheme.
We call a divisor $D$ \emph{coaffine},
if for every linearly equivalent effective divisor $E$
the complement of $E$ is affine.

$D$ is called \emph{strongly coaffine}, if $nD$ is coaffine or trivial
for every $n \in \ZZ$.
\end{definition}

\begin{lemma}
\label{coaffinefund}
Let $X$ be a normal scheme.
Then the following hold.

\renewcommand{\labelenumi}{(\roman{enumi})}
\begin{enumerate}
\item
The trivial divisor is coaffine if and only if $X$ is affine.
\item
If a positive multiple of $D$ is coaffine, then $D$ itself is coaffine.
\item
Suppose that $X$ is affine. Then $D$ is coaffine if and only if
$D_x \subseteq \Spec \, \O_x$ is coaffine for every point $x \in X$.
Cartier divisors are coaffine.
\item
If $W \subseteq X$ contains all points of codimension one and
if $D$ is coaffine on $W$, then $D$ is coaffine on $X$.
\end{enumerate}
\end{lemma}
\proof
(i) is clear.

(ii). Let $D+(q)$ be effective, $q \in K(X)$.
Then $n(D+(q))=nD+(q^n)$ is effective with the same support and its complement
is affine.
 
(iii). Let $X=\Spec \,A$.
If $D$ is not coaffine we find $E=D+(q)$ effective such that the complement
is not affine.
Since the affineness of an open subset $U \subseteq X$
in an affine scheme $X$ is a local property in $X$,
we find a point $x \in X$ where the complement of
$E_x$ is not affine, and
$E_x=(D+(q))_x=D_x+(q)_x$, thus $D_x$ is not coaffine.
 
For the converse, suppose that $D_x \subset \Spec \O_x$ is not coaffine, let
$D_x +(q)_x$ be effective in $\Spec \,\O_x$ such that the complement
is not affine.
Replacing $D$ by $E=D+(q)$ we find a divisor which is
locally in $x$ effective with non-affine complement.
Let $\primid$ be a prime divisor of $E$ not passing through $x$.
We find a global function $f \in A$, which is a unit in $x$ and
zero in $\primid$. Using such functions we obtain
an effective divisor equivalent with $E$ without
changing the localization at $x$.

(iv). Suppose that $D$ is effective. Then
$X-\Supp_X\, D \supseteq W-\Supp_W \,  D=U$.
If $U$ is affine, the other inclusion is also true,
since $W$ contains all points of codimension one.
Thus $X-\Supp_X\, D$ is affine.
\qed
\medskip

We come to the definition of the affine class group of a normal scheme.
Since the coaffine divisors do not form a subgroup of the divisor class group
in general,
our idea is to look at divisors which do not change
affineness properties of other divisors.

\begin{definition}
We say that a divisor $D$ is \emph{affine trivial}
if for every strongly coaffine divisor
$E$, the divisor $D+E$ is again strongly coaffine.
\end{definition}

\begin{proposition}
\label{afftrivfund}
Let $X$ be a normal scheme. Then the following hold.
\renewcommand{\labelenumi}{(\roman{enumi})}

\begin{enumerate}
\item
The set of affine trivial divisors form a subgroup
of $\Div X$. Every principal divisor is affine trivial.
\item
If $D$ is affine trivial, it is strongly coaffine.
\item
Suppose that $X$ is affine. If $kD$ is affine trivial, $k \neq 0$, then
$D$ is affine trivial.
\item
Suppose that $X$ is affine. If $D_x$ is affine trivial for every point
$x \in X$, then $D$ is affine trivial.
Every Cartier divisor is affine trivial.
\end{enumerate}
\end{proposition}
\proof
(i).
Suppose that $D$ and $D'$ are affine trivial, and let $E$ be strongly coaffine.
Then also $D+(D'+E)$ is strongly coaffine.
Let $D$ be affine trivial and $E$ strongly coaffine.
Then $-E$ is strongly coaffine and therefore
$D-E$ is strongly coaffine, thus $-D+E$ is strongly coaffine
and $-D$ is affine trivial.
The principal divisors are affine trivial, since coaffineness
is a property of the divisor class.

(ii) is clear.

(iii) Let $kD$ be affine trivial and let $E$ be strongly coaffine.
We may suppose $k >0$.
$kE$ is also strongly coaffine and therefore
$k(D+E)=kD+kE$ is strongly coaffine.
So $kn(D+E)$ is coaffine for $n \in \ZZ$ and thus
$n(D+E)$ is coaffine due to \ref{coaffinefund}.
 
(iv). 
In an affine scheme $X$ a divisor is (strongly) coaffine if and only if
it is locally (strongly) coaffine in every point.
Suppose $D_x$ is affine trivial for all $x \in X$ and let
$E$ be strongly coaffine.
Then $(D+E)_x=D_x+E_x$ is strongly coaffine and also $D+E$.
\qed

\begin{definition}
Let $X$ be a normal scheme.
We call the residue class group of the divisors modulo
the subgroup of affine trivial divisors the
\emph{affine class group} of $X$, denoted by
$$ ACl\, X \, .$$
\end{definition}

\begin{remark}
The affine class group is a residue class group of
the divisor class group. If $X$ is affine, it is a
torsion free factor group
of the Weil divisors modulo Cartier divisors.
\end{remark}

\begin{theorem}
Let $X$ be a normal scheme.
The affine class group of $X$ vanishes if and only if
every effective divisor is trivial or the
complement of its support is affine.
\par\smallskip\noindent
If $X$ is affine, then $ACl \, X=0$ if and only if the
complement of every hypersurface is affine.
\par\smallskip\noindent
If $X$ is the punctured scheme of a local ring or
a proper scheme over a field,
$ACl\, X=0$ holds if and only if the complement of every
non-empty hypersurface is affine.
\end{theorem}
\proof 
This is clear.
\qed

\begin{remark}
One important application of the property that the complement
of every hypersurface is affine is the theorem
of van der Waerden on the ramification of birational
morphisms, see \cite{EGAIV}, \S 21.12.
If $g:X \ra Y$ is a birational morphism of finite type and if
$Y$ is normal such that every point $y \in Y$ has this property
(or the property (W) of Grothendieck, \cite{EGAIV}, \S 21.12.8),
then the locus where
$g$ is not a local isomorphism has pure codimension one.
\end{remark}

%===========================================================
\section{Faithfully flat extensions}

We give a criterion for coaffine divisors.
The condition that $D$ is basepoint free in codimension one is always
fulfilled if $X$ is affine. $D$ corresponds to a reflexive module
$ U \longmapsto \shL_D(U)
=\{ q \in K : D+(q) \geq 0 \, \mbox{ on }\, U \} $,
which is an invertible sheaf outside the base locus.

\begin{lemma}
\label{coaffinecrit}
Let $X$ be a normal scheme, $A=\Gamma(X,\O_X)$
and let $D$ be a Weil divisor. 
Suppose that $D$ has no fixed components. Let $W$ denote
the complement of the base locus,
and let $\shL =\shL_D|_W$ be the corresponding invertible sheaf on $W$.
Then the following are equivalent.

\renewcommand{\labelenumi}{(\roman{enumi})}

\begin{enumerate}
\item
$D$ is coaffine on $X$.
\item
$D$ is coaffine on $W$.
\item
For every basepoint free linear subsystem
$(s_0,\ldots,s_n)$ of $\Gamma(W, \shL)$
the mapping $W \ra \PP_A^n$ is affine.
\item
There exists an $A$-generating system for the $A$-module
$\Gamma(X,\shL_D)=\Gamma(W,\shL)$
consisting of sections having affine complement.
\end{enumerate}
\end{lemma}
\proof
Let $E$ be effective and equivalent with $D$.
Then
$$X-\Supp_X\, E \supseteq (X-\Supp_X\, E)\cap W
= W-\Supp_W \,   E \, .$$
Since the base locus $X-W$ is contained in $\Supp_X\, E$,
both sets coincide. This gives (i) $\Rightarrow $ (ii),
the converse was proven in \ref{coaffinefund}.
 
(ii) $\Rightarrow $ (iii).
(ii) means that for every section
$s \in \shL_D(W)$ the complement of
the zero locus
$V(s)$ in $W$ is affine.
A basepoint free linear system $s_0,\ldots,s_n \in \Gamma(W, \shL)$
defines then an affine mapping $f: W \ra \PP^n_A$,
since $f^{-1}(D_+(x_i)) = W_{s_i}=W-V(s_i)$ is affine.
 
(iii) $\Rightarrow $ (iv) is clear. (iv) $\Rightarrow$ (ii).
Let $0 \neq s \in \shL(W)$ be a section corresponding
to the effective divisor $E$.
Then $s=a_0s_0+\ldots+a_ks_k$ where the $s_i$ are a generating system.
We may suppose that the linear system $(s_0,\ldots,s_k)$ is basepoint free.
Thus we get an affine mapping $\varphi: W \ra \PP_A^r$, showing that
$W_s= \varphi^{-1}(D_+(a_0x_0+\ldots+a_kx_k))$ is affine. 
\qed

\begin{remark}
The conditions in \ref{coaffinecrit} are not equivalent with
the property that there exists an affine morphism
$f: W \ra \PP_A^n$
with $\shL_D=f^\ast \O(1)$.
Consider for example the structure sheaf on a quasi affine, non-affine
scheme.
\end{remark}

\begin{lemma}
\label{coaffinepullback}
Let $p:X'=\Spec \,A' \ra X=\Spec \,A$
be flat, where $A$ and $A'$ are normal noetherian domains.
Let $D$ be coaffine on $\Spec \,A$.
Then the pull-back $D'=p^*D$ is also coaffine.
If the mapping is faithfully flat, the converse holds as well.
\end{lemma}
\proof
Let $X^{(1)}$ denote the set of prime ideals of height one
and let $W \subseteq X$ denote the basepoint free locus
of $D$. Thus $X^{(1)} \subseteq W$ and $\shL=\shL_D$ is invertible on $W$,
let $W'=p^{-1}(W)$.
We have $X^{'(1)} \subseteq W'$ since going down holds for flat mappings.
Due to \ref{coaffinefund} it suffice to show that
$\shL'$ is coaffine on $W'$.
 
Flatness gives $\Gamma(W',\shL')=\Gamma(W,\shL) \otimes_AA'$
and an $A$-generating system $s_i,\, i \in I,$ of
$\Gamma(W,\shL)$ gives an
$A'$-generating system
$s_i'=s_i \otimes 1,\, i \in I,$
of $\shL'(W')$.
The sections $s_i'$ have again affine complement
and the criterion \ref{coaffinecrit} (iv) shows that $D'$ is coaffine.
 
Suppose now that the mapping is faithfully flat
and that $p^*D$ is coaffine in $\Spec \,A'$.
Let $E$ be an effective representative of $D$.
Then $p^*E$ is effective with affine complement
and its support is the preimage of the support of
$E$. An open subset is affine if and only if its preimage under a
faithfully flat morphism is affine. Thus $X-\Supp\, E$ is affine.   
\qed

\begin{theorem}
\label{maintheorem}
Let $A \ra A' $ be a faithfully flat extension
of normal noetherian domains $A$ and $A'$.
Suppose that the mapping
$p^*: Cl\, A \ra Cl\, A'$ is surjective.
Then $ACl \, A= ACl \, A'$.
\end{theorem}
\proof
Let $D$ be affine trivial and let $E'$ be strongly coaffine.
We may write $E'=p^*E$ (as divisor class)
and $E$ is due to \ref{coaffinepullback} strongly coaffine.
Thus $p^*D+E'$ is strongly coaffine and $p^*D$ is affine trivial.
This gives a surjective morphism
$\, ACl\, A \ra ACl\, A'$.
A similiar argument shows that it is also injective.
\qed

\begin{corollary}
\label{affinclassgauss}
Let $A$ be a normal noetherian domain.
Then the affine class groups of $A$ and of the polynomial ring
$A[T]$ coincide.
\end{corollary}
\proof
The extension is faithfully flat and the divisor class groups
are the same.
\qed

\begin{corollary}
Let $A$ be a normal noetherian domain and let $E$ be a vector bundle
over $\Spec\, A$.
Then the affine class groups of $\, \Spec\, A$ and of $E$ coincide.
\end{corollary}
\proof
This is also clear.
\qed
\medskip

The following corollary answers one of the motivating questions of this paper.
It is analogous to the theorem of Gauss, that the polynomial ring
inherits factoriality from the base.

\begin{corollary}
\label{zylindernull}
Let $A$ be a normal noetherian domain such that
the complement of every hypersurface in $\Spec \, A$ is affine.
Then this is also true in $\Spec\, A[T]$.
\qed
\end{corollary}

\begin{remark}
If in $\Spec\, A$ every hypersurface is coaffine,
this is not true in the projective space $\PP^n_A$, because
the preimages of hypersurfaces have non-affine complement.
However, every hypersurface of $\PP^n_A$ which dominates the base
$\Spec \,A$ has an affine complement.
For this can be tested in the punctured affine cone and
the vertex of the cone lies in the closure of the preimage
of a dominating hypersurface.
\end{remark}

\begin{example}
\label{zylindercurve}
Cor. \ref{zylindernull} is not true without the condition that $A$ is normal.
Identify on the affine line $\AA_K^1$ two points and call the resulting
curve $C=\Spec\, A$. The complement of points on a one-dimensional
noetherian affine scheme is again affine.
Consider $S=C \times \AA_K^1=\Spec \,A[T]$.
$S$ is obtained by identifying two parallel lines on the affine plane.
The image curve of a skew line can not have an affine complement,
for its preimage consists of this line and two isolated points.
\end{example}

\begin{corollary}
\label{nagatazylinder}
Let $A$ be a two-dimensional normal excellent domain.
Then the complement of every {\rm(}hyper-{\rm)}surface in
$\Spec \,A[T]$ is affine.
The intersection of two surfaces in $\Spec\,A[T]$
has no isolated points.
\end{corollary}
\proof
This follows now from the Theorem of Nagata.
The intersection property follows from the general fact that
a hypersurface having affine complement has pure codimension
one on every closed subscheme
(see \cite{EGAIV}, \S 21.12. and \cite{brenner}).
\qed

\begin{example}
Cor. \ref{nagatazylinder} is again not true if $A$ is not normal.
Let $C$ be the curve of example \ref{zylindercurve}
and consider $C \times \AA_K ^2=S \times \AA^1_K$.
This three-dimensional variety arises by identifying
two parallel planes $P_1$ and $P_2$ in $\AA^3_K$.
One may find two disjoint surfaces $S_1$ and $S_2$ in affine space such that
$S_1 \cap P_1 =\emptyset , S_2 \cap P_2 =\emptyset$ and such that
the intersecting curves $S_1 \cap P_2$ and $S_2 \cap P_1$ have
only single points in common after identifying $P_1$ and $P_2$ .
The images of these two surfaces intersect then in isolated points on
$S \times \AA^1$.
\end{example}

\begin{problem}
The analogue of the theorem of Nagata in complex geometry states that
the complement of an analytic curve on a normal Stein surface $S$ is again
Stein (Theorem of Simha, \cite{simha}).
Analogue to \ref{nagatazylinder} one may ask: is the complement of every
analytic hypersurface in $S \times \CC$ Stein?
\end{problem}

We consider now the affine class group of a formal
power series ring $A[[T]]$.
A normal domain $A$ is said to have a
\emph{discrete divisor class group} if the
mapping $Cl\, A \ra Cl\, A[[T]]$ is
bijective, see \cite{fossum}, \S 19.
This holds for a normal excellent $k$-algebra $A$ over a field
of characterisic zero if and only if $A$ is 1-\emph{rational},
meaning that $H^1(X',\O_{X'})=0$, where $X' \ra \Spec \,A$
is a resolution of singularities, see \cite{binsto2}, 6.1.

\begin{corollary}
If the normal domain $A$ has discrete divisor class group,
the affine class groups of $A$ and of $A[[T]]$ coincide.
\qed
\end{corollary}

\begin{problem}
It is not clear whether $ACl\, A=0$ implies
$ACl\, A[[T]]=0 $.
This is not even clear if $A$ is factorial.
\end{problem}

\begin{example}
Even if a local domain $A$ is factorial, the affine class group
of the completion does not vanish in general.
Let
$$ A=K[U,V,W][X,Y]/(XY-(U^2-(1+W)V^2)) \, .$$
$A$ is factorial, because the polynomial $F=U^2-(1+W)V^2$ is irreducible
in $K[U,V,W]$. Over the completion $K[[U,V,W]]$ we find
$F=(U- \sqrt{1+W}V)(U+ \sqrt{1+W}V)$ and in the
completion of $A$,
$$ \hat{A}= K[[U,V,W,X,Y]]/(XY-(U- \sqrt{1+W}V)(U+ \sqrt{1+W}V))$$
the prime ideals
$( X, U-\sqrt{1+W}V )$ and $(Y,U +\sqrt{1+W}V )$ are of height one,
but their sum is $(X,Y,U,V)$, so they meet in the closed point,
thus they cannot have affine complement ($\Char \, K \neq 2$).
\end{example}

\begin{example}
\label{complexification}
In a similar way we may give examples of
factorial domains such that after changing the base field
the affine class group does not vanish anymore.
For example, the $\RR$-Algebra
$A=\RR [U,V,X,Y]/(XY-(U^2+V^2)) $
is factorial, but $ACl \, A\otimes_\RR \CC \neq 0$.
\end{example}

%===========================================================
\section{The affine class group of products}
\par
\bigskip
\noindent
Generalizing the situation
$X \times \AA \ra X$ of the previous section, we study
the affine class group of an affine normal variety $X$
over an algebraically closed field $k$
in relation to the affine class group of the product
with another affine normal variety, in particular with an affine smooth curve
$C$.
Of course, if $X$ and $Y$ are both smooth, then their product is again
smooth and the affine class group vanishes as well.
First we derive some corollaries from \ref{maintheorem}

\begin{corollary}
\label{fieldextension}
Let $X=\Spec\, A$ be an affine normal variety over an
algebraically closed field $k$ of characteristic zero,
and let $k \subseteq K$ be a field extension.
If $Cl\, A$ is finitely generated,
then the affine class group of $X$ and of $X_K=\Spec \, A \otimes_kK$
coincide.
\end{corollary}
\proof
Under these conditions the divisor class groups
coincide, see \cite{binsto2}, 15.7, so the result follows from
\ref{maintheorem}.
\qed

\begin{remark}
This is not true if $k$ is not algebraically closed, as example
\ref{complexification} has shown.
\end{remark}

\begin{corollary}
Let $X=\Spec \, A$ be an affine normal variety
over an algebraically closed field $k$ of characteristic zero,
and let $Y=\Spec\, B$ be a variety, where $B$ is factorial.
Then the affine class groups of $X$ and of $X \times Y$ coincide.
\end{corollary}
\proof
The divisor class groups coincide, see \cite{binsto2}, 15.10. 
\qed
\medskip

Let now $X=\Spec\, A$ be an affine normal variety and let $C$ be an affine
smooth curve over an algebraically closed field $k$.
The easiest hypersurfaces in $X \times C$ are graphs of morphisms
$X \ra C$ and more generally the closures of graphs
of rational mappings $X \dasharrow C$. We consider here also mappings to
smooth projective curves.

\begin{lemma}
\label{graphclosure}
Let $X$ be an affine normal variety and let
$C$ be an affine smooth curve or the projective line.
Let $\overline{C}$ be the smooth projective closure of $C$.
Let $f:X \dasharrow C$ be a rational mapping such that
$\overline{f}:X \supseteq U \ra \overline{C}$ is an affine
morphism, where $U$ is the locus where $\overline{f}$ is defined.
Then the closure of the graph of $f$ in $X\times C$ has an
affine complement.
\end{lemma}
\proof
Let $P_1,\ldots,P_n$ be the points in $\overline{C}$ not in $C$.
Let $\overline{\Gamma_f} \subseteq X \times C$ be the closure
of the graph of $f:X \dasharrow C$.
Then
$\overline{\Gamma_f}= \overline{\Gamma_{\overline{f}}} \, \cap \, X \times C$.
Thus our set equals

\begin{eqnarray*}
X \times C - \overline{\Gamma_f} &=&
X \times C - \overline{\Gamma_{\overline{f}}}\, \cap \, X \times C \\
& = & X \times \overline{C}
- \overline{\Gamma_{\overline{f}}} - X \times \{ P_1,\ldots,P_n \} \\
& =& X \times \overline{C} -
\overline{\Gamma_{\overline{f}} \cup U \times \{ P_1,\ldots,P_n \} } \,  .
\end{eqnarray*}
Since $\overline{f}$ is supposed to be affine, this is also true
for
$ \overline{f} \times id : U \times \overline{C} \ra
\overline{C} \times \overline{C}$.
The preimage of 
$\triangle \, \cup \, (\overline{C} \times \{P_1,\ldots,P_n \})$
under $\overline{f} \times id$ is
$\Gamma_{\overline{f}} \cup (U \times \{ P_1,\ldots,P_n \})$.
If $C=\PP$, the complement of the diagonal in
$\PP \times \PP$ is affine.
Otherwise $n \geq 1$ and then the curve
$\triangle \cup (\overline{C} \times \{ P_1,\ldots,P_n \})$
in $\overline{C}  \times \overline{C}$ is ample, so its complement is affine.
Therefore
$U \times \overline{C}
- (\Gamma_{\overline{f}} \cup (U \times \{ P_1,\ldots,P_n \}))$
is affine.
Due to the following lemma this set equals $ X \times \overline{C} -
\overline{\Gamma_{\overline{f}} \cup (U \times \{ P_1,\ldots,P_n \}) } $.
\qed

\begin{lemma}
Let $V$ and $Y$ be normal varieties and let
$\varphi : V \supseteq W \ra Y$ be a rational map.
Suppose that $W$ is the maximal locus where $\varphi$ is defined and
that $W$ contains all points of codimension one.
Suppose $Z \subseteq Y$ with $Y-Z$ affine.
Then $W - \varphi^{-1}(Z) = V- \overline{\varphi^{-1} (Z)}$.
\end{lemma}
\proof
Consider
$ V - \overline{\varphi^{-1}(Z)} \supseteq
W - \varphi^{-1}(Z) \ra Y-Z $.
$W - \varphi^{-1}(Z)$ contains all points of height one of
$V - \overline{\varphi^{-1}(Z)}$. Since the target is affine and $V$ normal,
this mapping is defined on $V - \overline{\varphi^{-1}(Z)}$.
\qed

\begin{remark}
If $X$ is an affine variety, $C$ any curve and
$f:X \ra C$ a morphism, then the
complement of the graph in $X \times C$ is always affine.
For this it suffice to show that the mapping
$X \times C -\Gamma_f \hookrightarrow X \times C \ra X$
is affine, and this is locally true.
\end{remark}

\begin{corollary}
\label{graphclosure2}
Let $X$ and $C$ be as in {\rm \ref{graphclosure}}
and suppose that $ACl \, X= 0$.
Then the  complement of the closure of the graph
of a rational mapping $f:X \dasharrow C$ is affine.
\end{corollary}
\proof
We have only to show that $\overline{f}$ is affine.
This is clear if $f$ is constant, otherwise
the preimage of a point is a hypersurface.
\qed

\begin{example}
The complement of the closure of the graph of a
rational mapping
to a smooth projective curve $C$ of higher genus
need not be affine, because then the
diagonal is not ample.
Let $C={\rm Proj}\, A$ be an elliptic curve
where $A$ is normal and let $f: X=\Spec \,A \dasharrow C$ be the cone mapping
defined on $U=D(A_+)$.
Then the complement of the graph of $f$ in $U \times C$ is not affine.
To show this, let $D$ be a curve on $C \times C$ disjoint to the diagonal,
say $D=(x,x+p)$, $p \neq 0$.
Inside the complement lies as a closed subscheme
the preimage of $D$, and this equals
$\{ (y,z) : z=f(y)+p \}$. This is
also the graph of a mapping on $U$ and therefore
isomorphic to $U$, which is not affine.
\end{example}

\begin{example}
Two affine normal varieties $X$ and $Y$
may have trivial affine class group,
but $ACl \, (X \times Y) \neq 0$.
For this consider the two-dimensional variety $X=\Spec \,A$
of the preceeding example.
Since $X \times C$ contains non-affine complements of hypersurfaces
(which dominate $X$),
this is also true for $X \times U$ and then also for $X \times X$.
If we replace one factor by a suitable affine open subset $Y$ in the
blowing-up of $X$ we find even examples $X \times Y$
where one factor is smooth.
\end{example}

\begin{theorem}
\label{surfacecurve}
Let $X=\Spec \,A$ be an affine normal surface and let
$C$ be an affine smooth curve both over an algebraically closed field $k$.
Then the complement of every hypersurface in $X \times C$
is affine, $ACl \, (X \times C)=0$.
\end{theorem}

\proof
We deduce this statement from \ref{graphclosure2} by considering
finite normal extensions $X' \ra X$.
Since $X'$ is again two-dimensional and normal, we have
$ACl\, X'=0$, so the complement of the closure of a graph in
$X' \times C$ is affine. Let $K=K(X)=Q(A)$.
 
Let $H=V(\primid) \subseteq X \times C$ be an irreducible hypersurface.
We may suppose that it dominates the base $X$.
Then the field extension of the generic points
$K \ra \kappa(\primid)=:L$ is finite.

Let $x_1,\ldots,x_n \in L$ be a $K$-algebra-generating system of $L$
with integral equations $F_1,\ldots,F_n \in K[T]$.
Due to the theorem of Kronecker, see \cite{scheja}, 54.10, there exists
a finite field extension $K \subseteq K'$, where all these polynomials $F_i$
split in linear factors.
Then every residue field of
$K' \otimes_KL=K'[X_1,\ldots,X_n]/(F_1,\ldots,F_n + \mbox{other relations} )$
is isomorphic to $K'$.
 
Let $A'$ be the normal closure of $A$ in $K'$ and $X'=\Spec \,A'$.
Since $A$ is excellent, it is a finite extension of $A$.
It suffice to show that the complement of the preimage of $H$
in $X' \times C$ is affine, due to the theorem of Chevalley, \cite{EGAII},
6.7.1.
Since going down is true for $X' \ra X$ and $X' \times C \ra X \times C$,
see \cite{nagloc}, I.10.13, the preimage
of $H=V(\primid)$ equals the closure of the fiber
of $\primid$.
The fiber over $\primid$ is $\Spec \,K' \otimes_KL$, its points
being $K'$-points.
So we have to show that
a hypersurface $H' \subseteq X' \times C $ with $K(H')=K(X')$ has
affine complement.
But then $H' \ra X'$ is generically an isomorphism and
there exists a rational mapping
$X' \dasharrow H' \hookrightarrow X' \times C \ra C$ such that
$H'$ is the closure of the graph of this mapping.
\qed

\begin{problem}
Suppose that $X$ is an affine normal variety and $C$ is an affine smooth
curve. Is it true that $ACl\, (X \times C)=ACl\, X$?
\end{problem}

%===========================================================
\section{Examples}
In this section we compute the affine class group of
hyperbolas over factorial domains,
of monoid rings and of determinantal rings.
In all these cases the result is that the affine class group
is just the divisor class group modulo torsion.
In all these examples we proove
that the complement of a certain divisor is not affine
by giving a closed subscheme
where the codimension of the intersection is bigger than one.

\medskip
{\it Hyperbolas}

\smallskip
\noindent
Let $R$ be a noetherian factorial domain,
let $U_1,\ldots,U_r$ be non associated prime elements of $R$.
We consider the hyperbola
$$A=R[X,Y]/(XY-f) \mbox{ where }
f=U_1^{d_1} \cdot \ldots \cdot U_r^{d_r},\, d_i > 0 \, . $$
${\primid}_i=(U_i,X)$ and $\foq_i =(U_i,Y)$ are prime ideals
of height one, the divisor class group of $A$ is
$Cl\, A =\ZZ^r/(d_1,\ldots,d_r)$,
generated by $\primid_i,\, 1 \leq i \leq r$.
In $Cl \, A$ the relations ${\primid}_i=-\foq_i$ hold.

\begin{proposition}
\label{aclhyperbolavanish}

For $A$ as above the following are equivalent.
\renewcommand{\labelenumi}{(\roman{enumi})}

\begin{enumerate}
\item
$ACl\, A=0$
\item
$D({\primid}_i)$ is affine for all $i=1,\ldots,r$.
\item
$U_i$ and $U_j$ generate the unit ideal in $R$ for $i \neq j$.
\item
$A$ is locally $\QQ$-factorial.
\end{enumerate}
\end{proposition}
\proof
(i) $\Rightarrow $ (ii) is clear.
(ii) $\Longrightarrow $ (iii).
Suppose that $U_i$ and $U_j,\, i \neq j,$ do not generate the unit ideal in $R$.
Then there exists in $A$
a prime ideal $\for \supseteq (X,Y,U_i,U_j)$,
and its height is at least three.
Therefore $V(\fop_i) \cap V(\foq_j)$ has codimension $\geq 2$
on $V(\foq_j)$ and $D(X,U_i)$ is not affine.

Suppose that (iii) holds.
Let $T_i$ denote the product $\prod_{j \neq i} U_j$.
Then $T_1,\ldots,T_r$ generate the unit ideal, thus $D(T_i)$ cover
$\Spec\, A$. On the other hand, $Cl \, A_{T_i}$ is finite, so (iv) holds.
(iv) $\Rightarrow$ (i) is clear.
\qed
\medskip

Suppose furtheron that $R$ is local.
If $r \geq 2$, then $ACl\, A \neq 0$.
We characterize now the divisor classes which are coaffine.

\begin{proposition}
\label{hyperbolacoaffine}
Let $R$ be a local noetherian factorial domain, $U_1,\ldots,U_r$
non associated prime elements in $R$,
$f=U_1^{d_1}\cdot \ldots \cdot U_r^{d_r}$ and
$A=R[X,Y]/(XY-f)$.
Then a Weil divisor $D$ of $A$ is coaffine if and only if it is
a principal divisor or equivalent with a divisor
$\, n_1{\primid}_1+\ldots+n_r{\primid}_r$, $0 < n_i < d_i$ for $i=1,\ldots,r$.
\end{proposition}
\proof
Every divisor of $A$ is equivalent with
$n_1{\primid}_1+\ldots+n_r{\primid}_r$.
Adding multiples of $d=(d_1,\ldots,d_r)$ we consider
the biggest representative $\leq d$.

Suppose first that at least one entry is not positive.
We order the index set such that
the first $k$ places are positive ($k \geq 1$),
the places
$k+1$ to $m$ are negative and the others are zero.
We replace the $\primid_i$ where $n_i < 0$ by $\foq_i$
and get an equivalent effective divisor with support
${\primid}_1,\ldots,{\primid}_k,\foq_{k+1},\ldots,\foq_m$, $1 \leq k \leq m$.

If $m=k$, the ideal
${\primid}_1 \cap \ldots \cap {\primid}_k + \foq_j, \, j > k,$ has
height $\geq 3$.
So suppose that at least one entry is also negative, $ m > k$.
Consider the ideal
$\for=(Y+U_1^{d_1}\cdot \ldots \cdot U_k^{d_k},
X+U_{k+1}^{d_{k+1}} \cdot \ldots \cdot U_r^{d_r})$.
This is a prime ideal of height one in $A$,
and
$\for
+{\primid}_1 \cap \ldots\cap {\primid}_k \cap
\foq_{k+1} \cap \ldots\cap \foq_m $
has again height $\geq 3$.
These divisors are not coaffine.

So suppose that all $n_i$ are positive.
If $(n_1,\ldots,n_r)=(d_1,\ldots,d_r)$, we have a coaffine
principal divisor.
If $n_i=d_i$ for at least one $i$, but not for all,
then considering $n-d$ and interchanging $\fop_i$ by $\foq_i$
yields a situation as before.
 
So suppose now that $0 < n < d$. We have to show that such a divisor
is coaffine. Note
$E= -n_1\foq_1-\ldots-n_r\foq_r \cong n_1{\primid}_1+\ldots+n_r{\primid}_r =D$
and consider the reflexive ideal
$$ \fob= \{ q \in Q(A) :\, E+ (q)\geq 0 \}
=\foq_1^{(n_1)} \cap \ldots \cap \foq_r^{(n_r)}\, . $$
Using the homogeneous mapping
$ R[X,U^dX^{-1}]= A \ra A_X=R[X,X^{-1}]$
we find that
$\fob= (U^n,Y) $.
The divisors
$E+(U^n)= n_1{\primid}_1+\ldots+n_r{\primid}_r$ and
$E+(Y)=d_1\foq_1+\ldots+d_r\foq_r-(n_1\foq_1+\ldots+n_r\foq_r)
=(d_1 - n_1)\foq_1 +\ldots+ (d_r-n_r) \foq_r$
have an affine complement, since $n_i$ and $d_i-n_i >0$ for all $i=1,\ldots,r$.
Then \ref{coaffinecrit} shows that $D$ is coaffine.
\qed

\begin{theorem}
Let $R$ and $A$ as in {\rm \ref{hyperbolacoaffine}}.
Let $D$ be a Weil divisor of $B$.
Then the following are equivalent.
\renewcommand{\labelenumi}{(\roman{enumi})}

\begin{enumerate}
\item
A multiple of $D$ is principal.
\item
$D$ is strongly coaffine.
\item
$D$ is affine trivial.
\end{enumerate}
\par
\smallskip
\noindent
The affine class group of $A$ is $Cl\, A$ modulo torsion.
\end{theorem}
\proof
We show (ii) $\Rightarrow $ (i).
Since $D$ is coaffine, $D$ is equivalent with $\sum_i n_i \fop_i$, where
$0 < (n_1,\ldots,n_r) < (d_1,\ldots,d_r)$.
We have to show that the fractions $d_i/n_i$ are all equal.
Suppose $d_1/n_1 \leq d_i/n_i$
for all $1 \leq i \leq r$.
Then
$d_1(n_1,\ldots,n_r) \leq n_1(d_1,\ldots,d_r)$, where equality holds in the
first place.
The divisor
$ n_1(d_1,\ldots,d_r) -d_1(n_1,\ldots,n_r)\geq 0$
is effective and equivalent with $-d_1D$.
Since it is coaffine and its first entry is zero,
it must be zero everywhere, hence it is a principal divisor.
\qed

\begin{example}
The sum of two coaffine divisors need not be coaffine.
Let $B=K[U_1,U_2][X,Y]/(XY-U_1^3U_2^3)$.
The divisor $(1,1)=\fop_1+\fop_2$ is affine trivial.
The divisor $(1,2)$ is coaffine, but $(1,1)+(1,2)=(2,3)=(-1,0)$ is not
coaffine.
Affine trivial divisors do not preserve
the property of being coaffine.
\end{example}

\medskip
{\it Monoid rings}

\smallskip
\noindent
We seek to compute the affine class group of a monoid ring $K[M]$,
where $K$ is a factorial domain and $M$ is a normal, finitely generated
torsion free monoid, see \cite{fulton}, \cite{bruns}.
Let $\Gamma=\Gamma(M)=\ZZ^d$ be
the quotient lattice and let $F_1,\ldots ,F_r$ be the facets of $M$.
Let $\fop_i=(T^m: m \not\in F_i)$ be the corresponding prime ideal of height
one and $\nu_i$ the corresponding valuation $\nu_i:\Gamma \ra \ZZ$.

The mapping $\nu=(\nu_1,\ldots,\nu_r): \Gamma \ra \ZZ^r$
is injective and $M=\Gamma \cap \NN^r$, hence it induces a canonical
mapping $K[M] \ra K[X_1,\ldots,X_r]$.
The factor group $\ZZ^r/\Gamma$ is the divisor class group of
$K[M]$. For an element $f \in M$ we have $f \not\in F_i$ iff $f \in \fop$
iff $\nu_i(f) >0$. We call
$\Supp\, (f)=\{F_i: f \not\in F_i \}$ the support of $f$.

\begin{lemma}
\label{monoidaffine}
Let $K[M]$ as above.
Then the complement of the effective divisor
$D=n_1{\primid}_1+\ldots+n_r{\primid}_r$ is affine if and only if
the support of $D$ is the support of a monomial $f \in M$.
\end{lemma}
\proof
Let $D=n_1F_1+\ldots+n_sF_s$,
$0 \leq s \leq r,\, n_i >0$
for $i=1,\ldots,s$ and
suppose that its support $\{F_1,\ldots,F_s \}$ is not the support of
a monomial function.
This means that for every monomial function
$f \in \foa ={\primid}_1 \cap \ldots \cap {\primid}_s$ there exists
another prime ideal
${\primid}_j,\, s<j \leq r,$ such that $f \in \fop_j$.
On the other hand there exists for every facet $F_i$ a monomial element $g_i$
such that $g_i \in F_i,\, g_i \not\in F_j$ for $j \neq i$.

We interpret these properties in terms of
the mapping $K[M] \subseteq K[X_1,\ldots,X_r]$
The extended ideal $\foa K[X_1,\ldots,X_r]$ is contained in
$(X_i)$ for $i=1,\ldots,s$.
Furthermore, it is contained in $(X_{s+1},\ldots,X_r)$.
But it cannot be contained in $(X_j)$, $s+1 \leq j \leq r$,
because for every $j$ there exists a monomial
$X_1^{\nu_1}\ldots X_r^{\nu_r}$ where $\nu_j=0$ and $\nu_i >0$ for
$i \neq j$. If $j \geq s+1$, this monomial belongs
to the extended ideal, but not to $(X_j)$.
So the extended ideal has not pure codimension one and its
complement cannot be affine.
\qed

\begin{theorem}
\label{monoidacl}
Let $D$ be a Weil divisor of the monoid ring $K[M]$, where $K$ is a factorial
domain. Then the following are equivalent.

\renewcommand{\labelenumi}{(\roman{enumi})}

\begin{enumerate}
\item
A multiple of $D$ is principal.
\item
$D$ is strongly coaffine.
\item
$D$ is affine trivial.
\end{enumerate}
\par
\smallskip
\noindent
The affine class group of $K[M]$ is $Cl\, K[M]$ modulo torsion.
\end{theorem}
\proof
Let
$D=n_1F_1+\ldots+n_rF_r$ be strongly coaffine, $n_i \geq 0$.
Due to \ref{monoidaffine}
the support of $D$ is the support of a function $f \in M$.
Thus there exist numbers $k,m \in \NN$ such that
$kD-m(f)$ is effective and such that the numbers of zeros has decreased.
This new divisor is equivalent with
$kD$ and therefore again coaffine.
Inductively we arrive at $t,s \in \NN$ such that $tD-s(g)=0$.
\qed

\begin{corollary}
Let $K[M]$ as before. Then the following are equivalent.
\renewcommand{\labelenumi}{(\roman{enumi})}

\begin{enumerate}
\item
$M$ is simplicial, meaning that the number of facets equals the
dimension.
\item
$K[M]$ is $\QQ$-factorial.
\item
$\Spec \,K[M]$ is a quotient singularity.
\item
$ACl \, K[M]=0$.
\end{enumerate}
\end{corollary}
\proof
(i) and (ii) are equivalent due to the explicit description
of $\Cl \, K[M]$ as the cokernel of
$\nu: \Gamma \ra \ZZ^r$.
Suppose (ii) holds.
$K[M] \ra K[X_1,\ldots,X_r]$ is the subring of degree zero
under the graduation given by $\ZZ^r \ra Cl \, K[M]$,
and this graduation
corresponds to the operation of the finite group scheme
$\Spec \,K[Cl \,(K[M])]$ on $\AA_K^r$, $\Spec \, K[M]$ being the
quotient.
(iii) $\Rightarrow $ (iv) follows from the
theorem of Chevalley \cite{EGAIV},6.7.1.
(iv) $\Rightarrow $ (ii) follows from \ref{monoidacl}.
\qed

\medskip
{\it Determinantal rings}

\smallskip
\noindent
Let $K$ denote a field and let $0<k\leq {\rm min}\,(m,n)$.
The $k$-minors of the $(m \times n)$-matrix
$$ \left( 
\begin{array}{ccc}
X_{11} & \ldots & X_{1n} \cr
\ldots & \ldots & \ldots \cr
X_{m1} & \ldots & X_{mn} \cr
\end{array}
\right) $$
define an ideal $I_k$ in the polynomial ring
$K[X_{ij}:\, 1 \leq i \leq m,\, 1 \leq j \leq n]$.

The ideals $I_k$ are prime of height $mn-(m+n-k+1)(k-1)$ and the
determinantal rings $R_k=K[X_{ij}]/I_k$ are normal
Cohen-Macaulay domains of dimension $(m+n-k+1)(k-1)$,
see \cite{bruns}, Theorem 7.3.1.

\begin{theorem}
Let $1< k \leq {\rm min}(m,n)$.
Let $\primid$ be the ideal of $R_k$ generated
by the $(k-1)$-minors of the first $(k-1)$ rows.
Then $\primid$ is a prime ideal of height one
and
$ACl \, R_k =Cl\, R_k =\ZZ $, generated by $\primid$.
\end{theorem}
\proof
For the second equality see \cite{bruns}, Theorem 7.3.5.
It suffice to show that
the complement of $V({\primid})$ is not affine.
Consider the mapping
$$K[X_{ij}:\, 1 \leq i \leq m,\, 1 \leq j \leq n]/I_k
\ra K[X_{ij}:\, 1 \leq i \leq k-1,\, 1 \leq j \leq n]$$
where $X_{ij} \mapsto X_{ij}$ for $i \leq k-1$ and
$X_{ij} \mapsto 0$ for $i \geq k$.
This is well defined.
The image
of $\primid$ is just the ideal $I_{k-1}$ of the
$(k-1)$-minors in the $(k-1)n$-dimensional polynomial ring.
Its height is
$(k-1)n-(k-1+n-(k-1)+1)(k-2)=nk-n-(n+1)(k-2)=n-k+2 \geq 2$,
thus $D(\primid)$ is not affine.
\qed

\begin{problem}
For which local normal domains $A$ does
the equation $ACl\, A =Cl\, A /{\rm torsion}$ hold?
Is this true for a rational singularity?
\end{problem}

%===========================================================
\section{Cones over projective varieties}

Let $K$ be a field and let
$A$ be an $\NN$-graded normal $K$-Algebra, finitely generated
by homogeneous elements of degree one.
Let $X={\rm Spec}\, A, \, X'=X-V(A_+)$ and $ Y={\rm Proj}\, A$
be the projective variety over $\Spec \, K$.
Let $H$ be the corresponding very ample divisor on $Y$.
The cone mapping
$p: X' \ra Y$ defines the pull-back $p^*: Div \, Y \ra Div \, X' =
Div \, X$.
We seek to relate the affine class group of the cone to the numerical
class group of the projective variety.

\begin{lemma}
Let $Y=\Proj \, A$ be as above, and $\dim \, Y \geq 1$.
Let $D$ be a
Weil divisor on $Y$.
Then the following hold.
\renewcommand{\labelenumi}{(\roman{enumi})}

\begin{enumerate}
\item
$p^*D$ is coaffine on $X$
if and only if for all
$k \in \ZZ$ the divisors $D+kH$ are coaffine on $Y$.
\item
$p^*D$ is strongly coaffine on $X$ if and only if
for all $m,n \in \ZZ$ the divisors
$nD+mH$ are coaffine or trivial.
\item
$p^*E$ is affine trivial on $X$ if and only
if for every divisor $D$ on $Y$ with the property that
for every $n,m \in \ZZ$ the divisors
$nD+mH$ are coaffine or trivial, this property holds for
$E+D$ as well.
\end{enumerate}
\end{lemma}
\proof
(i).
Let $\shL_D$ be the reflexive module on $Y$ corresponding to $D$.
Then
$\Gamma(X,p^*\shL)=\Gamma(X',p^* \shL)
= \bigoplus_{k \in \ZZ} \Gamma(Y,\shL_{D+kH})$.
Due to \ref{coaffinecrit},
$p^*D$ is coaffine on $X$ if for each homogeneous section
$s \in \Gamma(X,p^*\shL)_k=\Gamma(Y,\shL_{D+kH})$
the complement of the corresponding divisor on $X$ is affine, and this is
the case if and only if each effective divisor $D' \cong D+kH$ on $Y$
has an affine complement.

(ii) and (iii) follow from (i).
\qed

\begin{corollary}
Let $Y={\rm Proj}\, A$ be a projective normal variety over a field
$K$, where $A$ is a normal graded $K$-Algebra, finitely generated
by forms of degree one.
Then
$ACl \, A=0 $ if and only if the complement of every non-empty
hypersurface on $Y$ is affine.
\qed
\end{corollary}

\begin{proposition}
\label{numtrivial}
Let $Y=\Proj\, A$ be a smooth projective variety over an algebraically
closed field, $H$ the corresponding
very ample divisor.
Let $D$ be a Weil divisor on $Y$.
Suppose that $D$ and $H$ are linearly dependent
in the numerical class group
${\rm Num}\, Y$.
Then $p^*D$ is strongly coaffine.
\par
\smallskip
\noindent
Suppose furthermore that in $\Spec \, A$ every strongly coaffine divisor
is affine trivial.
Then
$ACl\, A$ is a factor group of
${\rm Num}\, Y$ and it is finitely generated.
\end{proposition}
\proof
Let $D \equiv qH$ in
${\rm Num}\, Y \otimes \QQ$, $q \in \QQ$.
We have to show that for $n,m \in \ZZ$ the divisors
$nD+mH$ are coaffine or trivial on $Y$.
We have $nD+mH \equiv rH,\, r \in \QQ,$and therefore
$p(nD+mH) \equiv kH,\, 0 \neq p \in \NN,\, k \in \ZZ$.
Suppose $k >0$. Ampelness is due to the criterion
of Seshadri (see \cite{kollar}, VI.2.18) a property of the numerical class,
thus $nD+mH$ is ample, hence coaffine.
If $k = 0$, the divisor $nD+mH$ is numerically trivial, hence
trivial or it has no sections $\neq 0$ at all.
If $k<0$, the divisor has no sections $\neq 0$.
\qed

\begin{problem}
Is the pull back of a numerically trivial divisor affine trivial?
Is the affine class group of a normal domain finitely generated?
\end{problem}

Let $\pi:Y \ra C$ be a geometrically ruled surface over a base curve $C$.
The Picard group of $Y$ is given by
$\Pic\, Y \cong \ZZ \oplus \pi^*({\rm Pic}\, C)$,
and the numerical class group is
${\rm Num}\, Y= \ZZ \oplus \ZZ$, \cite{haralg}, Prop. V.2.3.
We will compute the affine class group of an affine cone over $Y$.

\begin{proposition}
Let $Y$ be a geometric ruled surface over the curve $C$,
$A$ a normal homogeneous coordinate ring for $Y$ with
corresponding very ample sheaf $H$.
Then a divisor $p^*D$ is strongly coaffine in $\Spec\, A$ if and only if
$D$ is numerically dependent with $H$.
\par
\smallskip\noindent
The affine class group of the cone is $ACl\, A=\ZZ  $.
\end{proposition}
\proof
The if part was proven in \ref{numtrivial}.
So suppose
that $D$ and $H$ are numerically independent.
Then there exists a linear combination $mD+nH$ of
numerical type $(0,k)$, where $k$ is positive.
Such a divisor is the pull back of a divisor on the base curve
of positive degree $k$,
hence a multiple of it is linearly equivalent to an effective divisor.
The complement of an effective divisor coming from the base curve
contains projective lines, hence it can not be affine.

It follows that the strongly coaffine divisors of $A$ form a subgroup
in the divisor class group, thus they are affine trivial.
Hence
$ACl\, A =(\Num\, Y/H) /{\rm torsion} \cong \ZZ\, .$
\qed

%===========================================================

\end{document}